\documentclass[11pt,a4paper]{article}
\usepackage{latexsym}
\usepackage{epsfig}
\usepackage{here}
\usepackage{amsmath}
\usepackage{amssymb}
\usepackage{amsfonts}
\usepackage{amsbsy}
\usepackage{graphicx}
\usepackage{graphics}
\usepackage{euscript}
\usepackage{epsfig}
\usepackage{varioref}
\usepackage{enumerate}
\usepackage{verbatim}
\usepackage{here}
\usepackage{makeidx}
\usepackage{color}

\makeindex

\def\R{\mathbb{R}}
\def\C{\mathbb{C}}
\def\Z{\mathbb{Z}}

\def\N{\mathbb{N}}
\def\u{{\bf u}}
\def\v{{\bf v}}
\def\j{{\bf j}}
\def\k{{\bf k}}

\def\x{{\bf x}}

\newtheorem{definition}{Definition}[section]
\newtheorem{theorem}{Theorem}[section]

\newtheorem{example}{Example}[section]

\newtheorem{remark}{Remark}[section]

\numberwithin{equation}{section}

\makeatletter
\renewcommand\section{\@startsection {section}{1}{\z@}%
 				   {-2.5ex \@plus -1ex \@minus -.2ex}%
                                   {1.3ex \@plus.2ex}%
                                   {\normalfont\large\bfseries}}
\renewcommand\subsection{\@startsection {subsection}{1}{\z@}%
 				   {-2.ex \@plus -1ex \@minus -.2ex}%
                                   {0.5ex \@plus.2ex}%
                                   {\normalfont\normalsize\bfseries}}
\makeatother

\begin{document}

\title{Shannon wavelet approximations of linear differential operators}

\author{
Erwan Deriaz\thanks{Institute of Mathematics, Polish Academy of Science}
\\ {Erwan.Deriaz@impan.gov.pl}}

\bigskip
\bigskip

\date{\today}

\maketitle

\vspace{-7cm}
Preprint IMPAN

\vspace{1mm}

\hrule
\vspace{7cm}

\begin{abstract}

Recent works emphasized the interest of numerical solution of PDE's with wavelets. In their works
\cite{CDV01,CDV02}, A.~Cohen, W.~Dahmen and R.~DeVore focussed on the non linear approximation aspect
of the wavelet approximation of PDE's to prove the relevance of such methods. In order to extend these results,
we focuss on the convergence of the iterative algorithm, and we consider
different possibilities offered by the wavelet theory: the tensorial wavelets and the derivation/integration
of wavelet bases. We also investigate the use of wavelet packets. We apply these extended results to
prove in the case of the Shannon wavelets, the convergence of the algorithm introduced in \cite{DP06a}.
This algorithm carries out the Leray projector with divergence-free wavelets.

\end{abstract}

\section{Introduction}

Since the end of the 80's, the mathematical \emph{Theory of Wavelets} has
invented new tools for numerical simulations. Thanks to the
Fast Wavelet Transform, wavelets provide efficient algorithms
including optimal preconditionners for elliptic operators \cite{Jaf92}.
Recently, Cohen-Dahmen-De Vore's articles \cite{CDV01,CDV02} demonstrated the optimal complexity of wavelet algorithms
for the solution of elliptic problems.

These works enhanced the interest for these methods applied to the solution of partial differential equations.
Wavelet approach also resulted in non-linear approximations
\cite{Coh03}, and in denoising methods \cite{Don94}.
Wavelet methods offer the possibility to regulate both the accuracy in space and the accuracy in frequency.

In the following, we apply the Shannon decomposition to differential operators in order to investigate the convergence of wavelet
algorithms to solve partial differential equations. In the first part, we recall some operator theory basic elements; we
indicate how the Shannon wavelets and wavelet packets can be used to part the support of the Fourier transform of
a function.

In the following part, we resume the works \cite{CDV01,CDV02} on the wavelet approximation of
differential operators and state the correlated theorem of convergence for constant coefficient operators. In the
frame of Shannon wavelets,
we show that this convergence depends on the chosen wavelet decomposition (MRA or tensorial).
We extend this study to Shannon wavelet packets. After that, we introduce a result on the derivation
of biorthogonal wavelets due to P.~G.~Lemari\'e-Rieusset \cite{Lem92} that serves to construct new wavelet approximations
of differential operators.

In the last two parts, we present explicit examples that are implicated in the numerical solution of the Navier-Stokes
equations: the solution of the implicit Laplacian $(Id-\alpha\Delta)^{-1}$ and the Leray projector ${\mathbb P}$.


\section{Symbol of an operator}

The convergence of the wavelet methods involves the partition of the spectra
of the operator. This partition is provided by the wavelet decomposition. Hence we need the notion of
symbol as introduced by Lars H\"ormander in his book \emph{The Analysis of Partial
Differencial Operators} \cite{Hor85}.

\ 

We denote by $\partial_j$ the derivation along the variable $x_j$, and by $D_j$ the
derivation $-i\partial_j$. For $\alpha=(\alpha_1,\dots,\alpha_d)\in \N^d$, we write
$D^\alpha=D_1^{\alpha_1}\dots D_d^{\alpha_d}$.
Let $u$ denote a Shwartz function of $d$ real variables (i.e. $u$ is $C^\infty$ and fast decreasing:
$\forall N\in\N,~\exists A>0~/~\forall x \in\mathbb{R}^d,~|u(x)|<A/(1+|x|^2)^{N/2}$).
We denote by $<\cdot,\cdot>$ the scalar product either on vectors either in dual spaces.

Then $\widehat{u}$ stands for the Fourier transform of $u$, i.e.
$$
\widehat{u}(\xi)=\int_{\xi\in \mathbb{R}^d} e^{-i<x,\xi>} u(x)dx
$$
We also denote by ${\mathcal F}$ the isomorphism of $L^2(\R^d)$ given by
${\mathcal F}:u\mapsto (2\pi)^{-d/2}\widehat{u}$.

The inverse Fourier transform is done by:
\begin{equation}
\label{invFour}
u(x)=(2\pi)^{-d/2}\int_{\xi\in \mathbb{R}^d} e^{i<x,\xi>} {\mathcal F}u(\xi)d\xi
\end{equation}
When we derivate the relation (\ref{invFour}), it yields:
$$
D^\alpha u(x)=(2\pi)^{-d/2}\int_{\xi\in \mathbb{R}^d} e^{i<x,\xi>} \xi^\alpha {\mathcal F}u(\xi)d\xi
$$
Thus derivating $u$ by $D^\alpha$ consists in multiplying the Fourier transform of $u$ by
$\xi^\alpha=\xi_1^{\alpha_1}\dots \xi_d^{\alpha_d}$. The function $\xi\mapsto \xi^\alpha$ is called the symbol of $D^\alpha$.
More generally speaking, if $a(\xi)$ is a $C^\infty$ function slowly increasing (i.e. such that
$\exists N\in\N,~A>0~/~\forall\xi\in\mathbb{R}^d,~|a(\xi)|<A(1+|\xi|^2)^{N/2}$), $a(D)$ defines an operator of
symbol $a(\xi)$ acting on the class of the Schwartz functions ${\mathcal S}$ by
\begin{equation}
a(D) u(x)=(2\pi)^{-d/2}\int_{\xi\in \mathbb{R}^d} e^{i<x,\xi>} a(\xi) {\mathcal F}u(\xi)d\xi
\end{equation}
Let's now consider a differential operator $P$ of order $m$ with variable coefficients $a_\alpha$ in ${\mathcal S}$,
$P=\sum_{|\alpha|=0}^{m}a_{\alpha}(x)D^{\alpha}$. Then,
instead of using the formula:
$$
Pu=(2\pi)^{-d/2}\int_{\xi\in \mathbb{R}^d} e^{i<x,\xi>} {\mathcal F}(Pu)(\xi)d\xi
$$
where
$$
{\mathcal F}(Pu)(\xi)=(2\pi)^{d/2}\sum_{|\alpha|=0}^{m}{\mathcal F}a_\alpha\ast
\xi^\alpha {\mathcal F}u
$$
that is no more a multiplication but an integral operator on ${\mathcal F}u$,
we use the formula (for $x\in\mathbb{R}^d$):
$$
Pu(x)=(2\pi)^{-d/2}\int_{\xi\in \mathbb{R}^d} e^{i<x,\xi>}\left(
\sum_{|\alpha|=0}^{m}a_{\alpha}(x)\xi^{\alpha} \right)
{\mathcal F}(u)(\xi)d\xi
$$
that we write:
\begin{equation}
\label{sympseudodiff}
Pu(x)=(2\pi)^{-d/2}\int_{\xi\in \mathbb{R}^d} e^{i<x,\xi>} p(x,\xi)
{\mathcal F}(u)(\xi)d\xi
\end{equation}
introducing the ``symbol'' $p(x,\xi)$ of $P$
$$
p(x,\xi)=\sum_{|\alpha|=0}^{m}a_{\alpha}(x)\xi^{\alpha}
$$
The formula (\ref{sympseudodiff}) gives us the possibility to define the operators $p(x,D)$
of symbol $p(x,\xi)$ that are not polynomials in $\xi$. These operators are called
pseudo-differential. The functions $p$ must verify regularity and increase properties
of polynomial type (see \cite{Hor85}).

We'll remark that the function ${\mathcal F}(Pu)$ is no more the function $p(x,\xi)
{\mathcal F}u(\xi)$ which appears in (\ref{sympseudodiff}) since the latter depends on $x$.

\ 

The following definition of an elliptic operator is given in \cite{Hor85}:

\begin{definition}[elliptic operators]
From the symbol $p(x,\xi)=\sum_{|\alpha|=0}^{m}a_{\alpha}(x)\xi^\alpha$ we extract
the principal symbol $p_m(x,\xi)=\sum_{|\alpha|=m}a_{\alpha}(x)\xi^\alpha$. A differential
operator $P$ is said to be elliptic iff
$$
\forall x\in\mathbb{R}^d,~~\forall \xi\in\mathbb{R}^d\setminus\{(0,\ldots,0)\},
~~~p_m(x,\xi)\neq 0
$$
\end{definition}

\ 

In the following, we'll need differential operators applied to vector functions
$\mathbb{R}^d\to \mathbb{R}^{m}$.
We denote by ${\bf u}$ with bold caracter the (multi-variable) vector function $u$
of real variables when it has several components.
For $\u$ having several components, the symbol $p(x,\xi)$ is a matrix: $\forall x,\xi\in\R^d$, $M(x,\xi)\in
\mathbb{C}^{n\times m}$.
Let $A=(A_{ij})_{1\leq i\leq n,~1\leq j\leq m}$ be a differential operator
$A:\left( H^s(\mathbb{R}^d)\right)^m\to \left( H^r(\mathbb{R}^d)\right)^n$, with:
$$
A_{ij}=\sum_{\alpha}a_{ij,\alpha}(x)D^{\alpha}
$$
Its symbol is $M=(m_{ij})_{1\leq i\leq n,~1\leq j\leq m}$ with
$$
m_{ij}(x,\xi)=\sum_{\alpha}a_{ij,\alpha}(x)\xi^{\alpha}
$$
We apply the operator $A$ componentwise as follows:
$$
(A{\bf u})_i=\sum_{j=1}^m m_{ij}(x,D)u_j
$$
Therefore, the multidimensional symbol can be handled in much the same way as in 1-D.

\begin{remark}
As the operator is applied to real functional spaces, its symbol verifies the same relation
as the Fourier transform of real functions, that is:
$$
\forall i,j~~~~~~~~m_{ij}(x,-\xi)=m_{ij}(x,\xi)^*
$$
where $z^*$ denotes the complex conjugate of $z$.
\end{remark}


\section{Shannon wavelet decomposition}

A good reference for the definition of the Shannon wavelets is Mallat's academic
book \cite{M98}. We first briefly recall the
construction of these wavelets.

The {\bf biorthogonal wavelets} are based on
{\bf scale filters $m$ and $n$} that provide the low-pass filter and the
high-pass filter. As the scale
function $\varphi(\frac{\cdot}{2})$ and the wavelet $\psi(\frac{\cdot}{2})$
belong to $V_0=\,{span}\{ \varphi(\cdot-k),~k\in{\mathbb Z} \}$, there exist two
sequences $(a_k)$ and $(b_k)$ such that:
$$
\varphi(\frac{x}{2})=\sum_{k\in\mathbb{Z}}a_k \varphi(x-k)~~,~~~~~~
\psi(\frac{x}{2})=\sum_{k\in\mathbb{Z}}b_k \varphi(x-k)
$$
We thus get after a Fourier Transform:
$$
\hat{\varphi}(2\xi)=m(\xi)\hat{\varphi}(\xi)~~,~~~~~~
\hat{\psi}(2\xi)=n(\xi)\hat{\varphi}(\xi)
$$
with $m(\xi)=\frac{1}{2}\sum_{k\in\mathbb{Z}}a_k e^{-ik.\xi}$~, ~~$n(\xi)=\frac{1}{2}\sum_{k\in\mathbb{Z}}b_k e^{-ik.\xi}$\\
\noindent The scale function is inferred from the filter by\,:
$$
\hat{\varphi}(\xi)=\hat{\varphi}(0)\prod_{j=1}^{\infty}m(\frac{\xi}{2^j})
$$
A wavelet basis $\{\psi_{jk}\}_{j,k\in\Z}$ with $\psi_{jk}(x)=2^{j/2}\psi(2^jx-k)$
forms a Riesz basis of $L^2(\R)$.\\
Similarly, if we denote by $H^t$ the set of distribution
functions $f$ such that $(1+|\xi|^2)^{t/2}\widehat{f}\in L^2$, $\{2^{tj}\psi_{jk}\}_{j,k\in\Z}$ provides a Riesz
basis of the Hilbert space $H^t(\R)$ for $-t_1\leq t \leq t_2$.
Where, if $(\tilde \phi,\tilde \psi)$ is the dual wavelet basis associated to $(\phi,\psi)$,
$t_1\in \R$ is the maximal number such that $\tilde \psi \in H^{t_1}(\R)$ and $t_2\in \R$
the maximal number such that $\psi \in H^{t_2}(\R)$.\\
We denote by $\ell^2_{t}$ the space of sequences $(u_{jk})_{(j,k)\in\Z^2}$ with the norm
$\|(u_{jk})_{(j,k)\in(\Z^2)}\|_{\ell^2_{t}}^2=\sum_{(j,k)\in(Z^2)} 2^{2tj} |u_{jk}|^{2}$.
One caracteristic that plays an important role in operator approximation is the semi-orthogonality
coefficients $B_t,b_t>0$ such that $\forall (u_{jk})_{j,k\in\Z}\in\ell^2_{t}$,
$$
b_t\sum_{j\in\Z}\| \sum_{k\in\Z} u_{jk} \psi_{jk} \|_{H^t} \leq \| \sum_{j,k\in\Z} u_{jk} \psi_{jk} \|_{H^t}
\leq B_t\sum_{j\in\Z}\| \sum_{k\in\Z} u_{jk} \psi_{jk} \|_{H^t}
$$

\ 

{\bf Shannon wavelets} have this particularity to have perfect 
low-pass and high-pass filters\,:
$$
m(\xi)=\left\{ \begin{array}{l} 1~~~~{\rm if}~~~\xi\in[-\frac{\pi}{2}+2k\pi,\frac{\pi}{2}+2k\pi],~~~k\in{\mathbb Z}
\\ 0~~~~{\rm if}~~~\xi\in[\frac{\pi}{2}+2k\pi,\frac{3\pi}{2}+2k\pi],~~~k\in{\mathbb Z} \end{array} \right.
$$
$$
n(\xi)=\left\{ \begin{array}{l} -e^{-i\xi}~~~~{\rm if}~~~\xi\in[\frac{\pi}{2}+2k\pi,\frac{3\pi}{2}+2k\pi],~~~k\in{\mathbb Z}
\\ 0~~~~{\rm if}~~~\xi\in[-\frac{\pi}{2}+2k\pi,\frac{\pi}{2}+2k\pi],~~~k\in{\mathbb Z} \end{array} \right.
$$
\begin{figure}
\begin{center}
\input{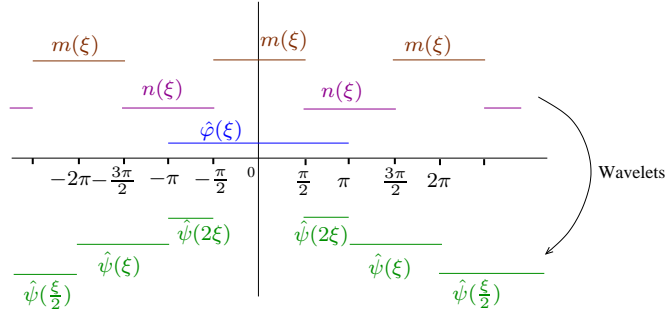}
\end{center}
\caption{Representation of the compact support of the Shannon filters and of the Fourier transform of the wavelets.
The vertical axis has no particular meaning, but the bars represent the compact support of the functions.}
\end{figure}
Then the corresponding scaling function writes:
$$
\hat\varphi(\xi)=\chi_{[-\pi,\pi]}(\xi)\qquad , \qquad \varphi(x)=\frac{\sin \pi x}{\pi x}
$$
and the wavelet:
$$
\hat\psi(\xi)=e^{-i\xi/2}~\chi_{[-2\pi,-\pi]\cup[\pi,2\pi]}(\xi)\qquad , \qquad
\psi(x)=\frac{\sin 2\pi(x-1/2)}{2\pi(x-1/2)}-\frac{\sin \pi(x-1/2)}{\pi(x-1/2)}
$$
where $\chi$ stands for the characteristic function i.e. $\chi_{E}(x)=1  ~{\rm if}~  x\in E,
0  ~{\rm if}~ x\notin E$.

\ 

In the multidimensional case, the {\bf tensorial Shannon decomposition}
can be written as follows:\\
Let ${\bf u}:\mathbb{R}^d\to \mathbb{R}^m$. The Shannon decomposition of $\u$ is given by:
\begin{equation}
\label{shaexpansion}
{\bf u}=\sum_{{\bf j}\in {\mathbb Z}^d}{\bf u}_{\bf j}
\end{equation}
with
$$
{\rm supp~}\widehat{{\bf u}_{\bf j}}\subset
\prod_{i=1}^d [-2^{j_i+1}\pi,-2^{j_i}\pi]\cup [2^{j_i}\pi,2^{j_i+1}\pi]
$$
For each scale parameter ${\bf j}\in {\mathbb Z}^d$, and for each component $\ell=1\dots m$, we have:
$$
u_{\ell,{\bf j}}(\x)=\sum_{{\bf k}\in{\mathbb Z}^d} 2^{\j/2}u_{\ell\,\j\k}\psi_1^\ell(2^{j_1}x_1-k_1)\dots
\psi_i^\ell(2^{j_i}x_{i}-k_{i})\dots\psi_d^\ell(2^{j_d}x_d-k_d)
$$
where $|\j|=\sum_i j_i$ and $\psi_i^\ell$ are wavelets of Shannon type, i.e.
${\rm supp\,}\widehat{\psi}_i^\ell\subset [-2\pi,-\pi]\cup[\pi,2\pi]$.

\subsection{Shannon wavelet packets}

With the above filters $m(\xi)$ and $n(\xi)$, we can define the Shannon wavelet packets.
The wavelet packets are defined by applying the filters $m$ and $n$ to the wavelets.
Hence we obtain two new wavelets $\psi_{(11)}$ and $\psi_{(10)}$ that are twice better
localised in the Fourier space (i.e. the compact supports of their Fourier
transforms are twice smaller):
\begin{equation}
\label{raf(10)}
\widehat{\psi_{(10)}}(2\xi)=m(\xi)\hat{\psi}(\xi)
\end{equation}
\begin{equation}
\label{raf(11)}
\widehat{\psi_{(11)}}(2\xi)=n(\xi)\hat{\psi}(\xi)
\end{equation}

The two of them are necessary to expend $L^2(\R)$, i.e. $L^2(\R)=\,{span}\{\psi_{(10))}(2^jx-2k),\psi_{(11)}(2^jx-2k)\}_{j,k\in\Z}$.
More precisely, the wavelet space at level $j$, $W_j$, admits 
$\{\psi_{(10))}(2^jx-2k),\psi_{(11)}(2^jx-2k) \}_{k\in\Z}$ as a Riesz basis.

\ 

\begin{figure}[!hb]
\begin{center}
\input{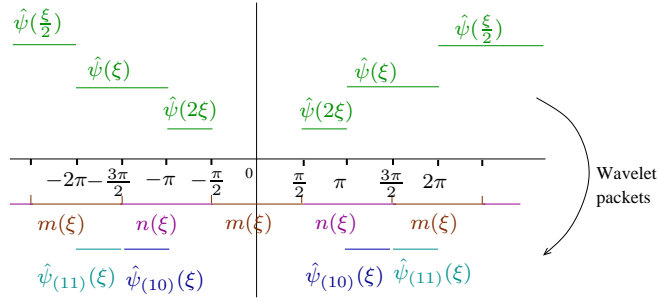}
\end{center}
\caption{The construction of Shannon wavelet packets represented in the Fourier space}
\end{figure}

The operations (\ref{raf(10)}) and (\ref{raf(11)}) on the wavelets can be iterated
as many times as one wants, and the Fourier support can be
shrunk as desired. In practice this operation can also be applied to usual wavelets but
doesn't come out with good results. Getting a better frequency localisation for usual
wavelet packets is still a challenging problem.

\ 

The Shannon wavelet packet decomposition gives us more oportunities to approximate any operators
than the classical Shannon decomposition
(\ref{shaexpansion}). We can report to part \ref{contwavpak} for their actual use.


\section{Solution of elliptic PDE's with wavelets -- the Richardson iteration}
\label{CDVtens}

In their paper \cite{CDV01}, A.~Cohen, W.~Dahmen and R.~DeVore consider a simple method to solve elliptic operator
equations with wavelets.
Let first consider $m=n=1$. In order to find the solution of the differential equation:
\begin{equation}
\label{eqnb}
Au=v
\end{equation}
where $A$ is a linear differential operator
and $u$ the unknown function, they use the expansions of $u$ and $v$ in wavelet bases. We denote by
$\underline{u}=(u_{jk})_{jk}$ the vector of wavelet coefficients:
$u_{jk}=<u,\tilde\psi_{jk}>$
with $\{\psi_{jk}\}_{j,k}$ and $\{\tilde\psi_{jk}\}_{j,k}$ two dual wavelet bases.
Then the expansion of $u$ writes $u=\sum_{j,k} u_{jk}\psi_{jk}$.\\

Let $\underline{A}$ be the variational discretisation of $A$ expressed in the wavelet basis $\{\psi_{jk}\}_{(j,k)}$
(it is called the Petrov-Galerkin stiffness matrix)
$\underline{A}=(<A\psi_{jk},\tilde\psi_{j'k'}>)_{j,j',k,k'}$ and $\underline{D}$
a wavelet preconditionner associated to the wavelet expansions (usually, it is the diagonal of $\underline{A}$, that
has the form Diag($2^{tj}$)). We assume that $A$ is continuous from $H^{t/2}$ to $H^{-t/2}$ and $a=<A\,\cdot,\cdot>$
is coercive.\\
In order to diversify the considered wavelet transforms in dimension $d>1$, let
$J$ denote the scale indice set ($J$ countable, e.g. $\Z\times\{0,1\}^{d*}$ for the MRA where $*$ means deprivated
of the element $(0,\dots,0)$, and $J=\Z^d$ for tensorial
wavelets). For a function $\alpha:J\to\R$ we define the $\ell^2_{\alpha}$ norm on the wavelet coefficients
$\underline{u}=(u_{\j \k})_{\j\in J,\, \k\in\Z^d}$ by
$$
\|(u_{\j \k})_{\j\in J,\, \k\in\Z^d}\|_{\ell^2_{\alpha}}=\sum_{\j\in J,\, \k\in\Z^d} 2^{\alpha(\j)} |u_{\j \k}|^{2}
$$
Usually, the $\ell^2_{t}$ corresponds to either the case $\alpha(\j)=t j_1$ in the MRA case $\j=(j_1,\varepsilon)\in
\Z\times\{0,1\}^{d*}$, either $\alpha(\j)=t \max(j_1,\dots,j_d)$ for the tensorial wavelets.\\

Then $\underline{D}^{-1}$ is continuous from $\ell^2_{-t/2}$ to $\ell^2_{t/2}$ (see \cite{Jaf92}).
We write the sequence $(\underline{u}_n)$, thanks to a Richardson iteration associated with a multiscale preconditionnning,
starting from $\underline{u}_0=0$\,:
\begin{equation}
\label{sequ}
\underline{u}_n=\underline{u}_{n-1}+\underline{D}^{-1}(\underline{v}-\underline{A}\,\underline{u}_{n-1})
\end{equation}
Then this method is said to converge if
$$
\exists \rho<1,~~~~\|\underline{v}-\underline{A}\,\underline{u}_n \|_{\ell^2_{-t/2}}\leq
\rho\|\underline{v}-\underline{A}\underline{u}_{n-1} \|_{\ell^2_{-t/2}}
$$
As we have:
$$
\underline{v}-\underline{A}\,\underline{u}_n=(\underline{Id}-\underline{A}\,\underline{D}^{-1})(\underline{v}-\underline{A}\,\underline{u}_{n-1})
$$
the algorithm converges if $\rho=\| \underline{Id}-\underline{A}\,\underline{D}^{-1}\| <1$, in the operator norm. And
$$
\underline{u}_n\to_{n\to\infty}\underline{u} \quad {\rm with} \quad \underline{A}\,\underline{u}=\underline{v}
$$
That is, as $u=\sum_{jk}u_{jk}\psi_{jk}\in H^{t/2}$,
$$
Au=v
$$

\ 

From  now on we think of $A$ as being an operator with constant coefficients. Let $m$, $n$ be unspecified natural numbers.
Hence we switch to vector spaces.
If we denote by $M(\xi)$ the symbol associated to $A$, we can express $A$
after a Fourier transform of the equation (\ref{eqnb}) as
$$
M(\xi) \widehat{\bf u}=\widehat{\bf v} ~~~~~~~~~~~\textrm{and the pseudo-inverse solution}~~~~~~~~~~~ \widehat{\bf u}
=M(\xi)^{\dag} \widehat{\bf v}
$$
with ${\bf u}\in(H^{t/2}({\mathbb R}^d))^m$, ${\bf v}\in(H^{-t/2}({\mathbb R}^d))^n$,
$M(\xi)\in{\mathbb C}^{n\times m}$ and $M(\xi)^{\dag}$ the pseudo-inverse of $M(\xi)$.\\
Remark that if $m=n$ and $M(\xi)M(\xi)^{\dag}=Id$, $M(\xi)^{\dag}=M(\xi)^{-1}$.

\ 

The idea for solving $A\u=\v$ is the following: we decompose ${\bf v}$ in a wavelet basis that splits the support of $\widehat{\bf v}$
$$
{\bf v}=\sum_{{\bf j}\in J}{\bf v}_{\bf j}
$$
If we denote by $({\bf e}_1,\dots,{\bf e}_n)$ the canonical basis of $\R^n$ then
${\bf v}_{\bf j}$ is the projection of $\v$ in the wavelet level
$W_\j=\,{span}(\{ {\bf \Psi}_{1,\j \k} \}_{\k\in\Z^d},\dots,\{ {\bf \Psi}_{n,\j \k} \}_{\k\in\Z^d})$
with each component $v_\ell {\bf e}_\ell$ of ${\bf v}$
decomposed in the wavelet basis $\{ {\bf \Psi}_{\ell,\j \k} \}_{\j\in J,\k\in\Z^d}$
(further we'll need this generalisation of \cite{CDV01} which uses an MRA).\\
For example in the tensorial case, we have $J=\Z^d$ and:
$$
v_{\ell \bf j}=\sum_{\k\in\Z^d}v_{\ell\j\k}\psi_{j_1 k_1}(x_1)\dots\psi_{j_d k_d}(x_d)
$$
This modification of the MRA case will prove usefull in part \ref{leraynum}.\\
Let us assume that for each ${\bf j}\in J$, $\widehat{\bf v}_{\bf j}$
and $\widehat{\bf u}_{\bf j}$ are compactly
supported (wavelet decompositions give us the
opportunity to do this with the desired accuracy). For each ${\bf j}\in J$, we build a
matrix $M_{\omega_{\bf j}}$ ($M_{\omega_{\bf j}}\in{\mathbb R}^{n\times m}$ depending on the compact support of
$\widehat{\bf v}_{\bf j}$) such that:
\begin{equation}
\label{approxop}
M_{\omega_{\bf j}}\thickapprox M(\xi)~~~~~~\textrm{for}~~\xi\in\textrm{supp}(\widehat{\bf v}_{\bf j})
\end{equation}
Then we approximate the relation $M(\xi) \widehat{\bf u}=\widehat{\bf v}$ by:
\begin{equation}
\label{appMxi}
\forall \j\in J,~\forall\k\in\Z^d,~~M_{\omega_{\bf j}}\left[
\begin{array}{c} u_{1\j\k} \\ \vdots \\ u_{n\j\k} \end{array} \right]
=\left[ \begin{array}{c} v_{1\j\k} \\ \vdots \\ v_{m\j\k} \end{array} \right]
\end{equation}

In the view of Richardson iteration, we take as a preconditionnner
$D=\sum_{\j\in J}M_{\omega_\j}P_\j$, where $P_\j$ is the projector on the wavelet level $W_\j$.
Then the corresponding discrete preconditionner $\underline{D}$ which applies to wavelet coefficients
is $\underline{D}=\sum_{\j\in J}M_{\omega_\j}\underline{P}_\j$
where $\underline{P}_\j$ is a diagonal matrix with ones on the lines $(\k,\j)_{\k\in\Z^d}$ and
zeros everywhere else.
In the case of tensorial wavelet basis ($J=\Z^d$), the space $W_\j$ is the $L^2$ closure of the
space generated by the family
$$
\{ {\bf \Psi}_{\ell,\j \k} \}_{1\leq \ell \leq m,\,\k\in\Z^d}=\{(\psi_{j_1 k_1}\times\dots\times\psi_{j_d k_d},
\underbrace{0,\dots,0}_{m-1})\}_{\k\in\Z^d}\cup
\dots\cup\{(\underbrace{0,\dots,0}_{m-1},\psi_{j_1 k_1}\times\dots\times\psi_{j_d k_d})\}_{\k\in\Z^d}
$$
In the following, we use the notation $M_{\omega}=\sum_{\j}M_{\omega_\j}P_\j$.\\
If we write the sequence (\ref{sequ}) with
$\underline{\v}_n=\underline{\v}-\underline{A}\,\underline{\u}_{n}$, it comes:
\begin{equation}
\label{sequw}
\underline{\u}_{0}=0~,~~~
\underline{\v}_{0}=\underline{\v}~,~~~~~~\underline{\u}_{n+1}=\underline{\u}_n+M_{\omega}^{\dag}\underline{\v}_{n}~~~~{\rm and}~~~~~
\underline{\v}_{n+1}=\underline{\v}_{n}-\underline{A}(\underline{\u}_{n+1}-\underline{\u}_{n})
\end{equation}

\begin{theorem}
\label{th1}
Let $M(\xi)$ be the symbol matrix associated to $A:H^{t/2}\to H^{-t/2}$ continuous. If the wavelet basis
$\{ {\bf \Psi}_{\ell,\j \k} \}_{1\leq \ell \leq m,\,\j\in J,\,\k\in\Z^d}$ provides a Riesz basis
of $H^{\pm t/2}$ (i.e. the associated decompositions $\v\mapsto \underline{\v},~H^{\pm t/2}\to \ell^2_{\pm t/2}$,
and reconstructions $\underline{\v}\mapsto \v,~\ell^2_{\pm t/2}\to H^{\pm t/2}$ are continuous).\\
Moreover, we suppose we have constructed for all $\j\in J$ matrices $M_{\omega_{\bf j}}\in {\mathbb R}^{n\times m}$
such that $M_{\omega}^\dag =\sum_{\j}M_{\omega_\j}^\dag P_\j:H^{-t/2}\to H^{t/2}$ is continuous. We also
assume that the wavelet decomposition $\v\mapsto ({\v}_\j)_{\j\in J}$ satisfies:
$$
\exists \tilde B>0~~~\textrm{such that}~~\forall{\v}\in H^{-t/2},~~~~~
\|(Id-A\, M_{\omega_{\bf j}}^\dag){\v}\|_{H^{-t/2}}^2\leq
\tilde B \sum_{\j\in J}\|(Id-A\, M_{\omega_{\bf j}}^\dag){\v}_\j\|_{H^{-t/2}}^2
$$
If there exist a real number $\rho\geq 0$ such that:
$$
\forall\j\in J,~~|||(Id-A\, M_{\omega_{\bf j}}^\dag)_{|W_\j}|||\leq \rho
$$
i.e.
$$
\forall\j\in J,~~~~\forall{\v}_\j\in W_\j,~~~~~
\|(Id-A\, M_{\omega_{\bf j}}^\dag){\v}_\j\|_{H^{-t/2}}\leq \rho\|{\v}_\j\|_{H^{-t/2}}
$$
then for $\rho$ small enough, the sequence $(\underline{\u}_n)_{n\in \N}$
defined by (\ref{sequw}) converges in $\ell^2_{t/2}$ to the wavelet expansion $\underline{\u}$ of $\u$ such that:
$$
\underline{A}\,\underline{\u}=\underline{\v}\quad  {\rm and~we~have} \quad A\u=\v
$$
\end{theorem}
\emph{proof:}\\
The graph of continuous operators can be summarized as follows:
$$
\begin{array}{ccc}
 & {\rm wavelet~transform} &  \\
\u\in {H^{t/2}} & \longleftrightarrow & \underline{\u}\in \ell^2_{t/2} \\
\ & & \\
A\downarrow ~~ \uparrow M_{\omega}^\dag &  & \underline{M}_{\omega}^\dag \uparrow~~ \downarrow \underline{A} \\
\ & & \\
\v=A\,\u\in{H^{-t/2}} & \longleftrightarrow & \underline{\v}=\underline{A}\,\underline{\u}\in\ell^2_{-t/2} \\
 & {\rm wavelet~transform} & 
\end{array}
$$
The operator $ M_{\omega}^\dag$ is not the inverse of $A$ but its approximation.\\

As the wavelet decompositions are continuous,
$$
\exists b,B>0~~~\textrm{such that}~~\forall{\v}\in H^{-t/2},~~~~~
b \sum_{\j\in J}\|{\v}_\j\|_{H^{-t/2}}^2\leq
\|{\v}\|_{H^{-t/2}}^2\leq
B \sum_{\j\in J}\|{\v}_\j\|_{H^{-t/2}}^2
$$
When $b=B=1$, the wavelet basis is said to be semi-orthogonal.\\
Then we have:
$$
\|\v_{n+1}\|_{H^{-t/2}}^2\leq \tilde B \sum_{\j\in J}\|(Id-A\, M_{\omega}^\dag){\v}_{n\,\j}\|_{H^{-t/2}}^2
\leq \tilde B \sum_{\j\in J}\rho^2 \|{\v}_{n\,\j}\|_{H^{-t/2}}^2
\leq \rho^2\frac{\tilde B}{b} \|{\v}_{n}\|_{H^{-t/2}}^2
$$

If $\rho^2 \tilde B/b<1$, as
$M_{\omega}^\dag$ is continuous, the serie
$\sum_n M_{\omega}^\dag \v_{n}$ converges in the Banach space $H^{t/2}$ to a solution $\u$ of the equation
$A\u=\v$.

\ 

The ideal wavelets that provide a minimal compact support for the Fourier transform are the
Shannon wavelets. In this case, as the compact supports of the Fourier transforms of wavelets from different levels are
disjoint, we have $\tilde B=b=1$ for all $A$. Shannon wavelets have an infinite support and are not used in practice.
But, in first approximation,
all wavelets behave as Shannon wavelets with more or less accuracy.
\begin{remark}
In the case of Shannon wavelets, as
$$
\widehat{\v}_{|\cup_\ell\,{\rm supp\,}(\widehat{{\bf \Psi}_{\ell\j\k}})}=\sum_{\ell,\k\in\Z^d} v_{\ell\j\k} \widehat{\bf \Psi}_{\ell\j\k}
$$
the equation \ref{appMxi} is equivalent to
$$
\forall \j\in J,~~~M_{\omega_\j}\widehat{\u}(\xi)=\widehat{\v}(\xi)~~~{\rm for}~~~ \xi\in \cup_\ell\,{\rm supp\,}(\widehat{{\bf \Psi}_{\ell\j\k}})
$$
In the futur, that will allow us to express this relation using the components $\v_\j$ of the Shannon
decomposition of $\v$.\\
One can also remark that as $M_{\omega}$ doesn't depend on $\xi$, we can apply the operator $M_{\omega}^\dag$ in the
physical space (expressed with wavelets) and not in the Fourier space.
\end{remark}


\section{Multiresolution analysis (MRA) versus tensorial basis}

\begin{figure}[!h]
\begin{center}
\input{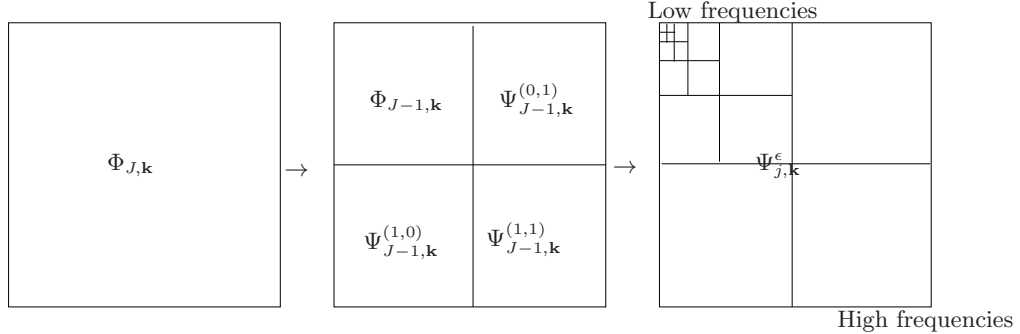}
\end{center}
\caption{\label{isotropic} Splitting of the Fourier modes induced by the 2D-MRA wavelet decomposition}
\end{figure}

\begin{figure}[!h]
\begin{center}
\input{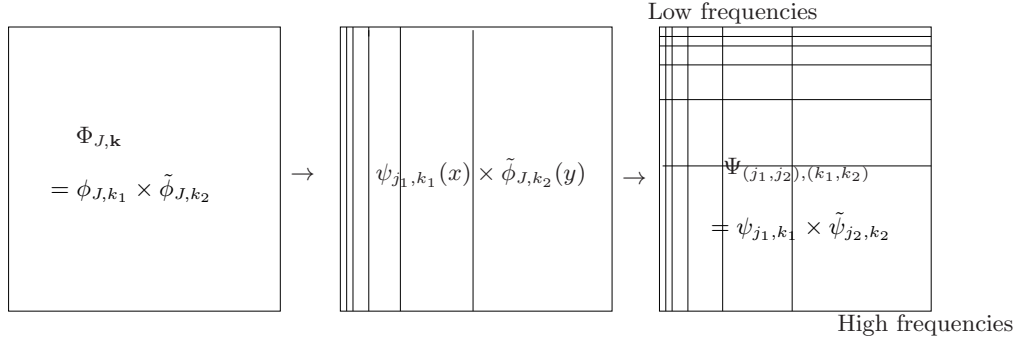} 
\end{center}
\caption{\label{anisotropic} Fourier splitting induced by the tensorial wavelet decomposition}
\end{figure}

There are two main different kinds of wavelet decompositions for a function on $\mathbb{R}^d$ with $d\geq 2$. It can be
decomposed either in a multidimensional multiresolution analysis or in a tensorial basis.
In an MRA, the wavelet decomposition of a function $u$ in 2D writes:
\begin{eqnarray*}
u(x_1,x_2)&=&\sum_{j\in {\mathbb Z}}\left( \sum_{(k_1,k_2)\in{\mathbb Z}^2}u_{j,k_1,k_2}^{(1,0)}~\psi_{0\, j,k_1}(x_1)~\varphi_{1\, j,k_2}(x_2)+
\sum_{(k_1,k_2)\in{\mathbb Z}^2}u_{j,k_1,k_2}^{(0,1)}~\varphi_{0\, j,k_1}(x_1)~\psi_{1\, j,k_2}(x_2)\right.\\
& & \left.\quad +\sum_{(k_1,k_2)\in{\mathbb Z}^2}u_{j,k_1,k_2}^{(1,1)}~\psi_{0\, j,k_1}(x_1)~\psi_{1\, j,k_2}(x_2)\right)
\end{eqnarray*}
where we used the notation $\psi_{j,k}(x)=\psi(2^j\,x-k)$.\\
While, in a tensorial basis it writes:
\begin{eqnarray*}
u(x_1,x_2) & =&\sum_{j_1\in{\mathbb Z}}\sum_{j_2\in{\mathbb Z}}\sum_{(k_1,k_2)\in{\mathbb Z}^2}
u_{j_1,j_2,k_1,k_2}~\psi_0(2^{j_1}x_1-k_1)~\psi_1(2^{j_2}x_2-k_2)
\end{eqnarray*}
These two decompositions correspond to two different partitions of the Fourier space (i.e. frequency domain). Both of them
are represented in figures \ref{isotropic} and \ref{anisotropic}. On each figure, in the last square, which corresponds
to the wavelet transform, the low frequencies are localised in the upper left corner of the square of coefficients,
and the high frequencies in the bottom right.

\subsection{Convergence theorems with Shannon wavelets}

To begin with, we'll only consider approximation matrices that are constant over each frequency domain indexed by
${\bf j}\in\Z^d$. 

The two previous decompositions induce different conditions for the approximation of the matrix $M(\xi)$.
Adding to part \ref{CDVtens}, we distinguish MRA and tensorial wavelet convergence theorems as follows:

\begin{theorem}[MRA]
If the symbol matix $M(\xi)$ admits a pseudo-inverse $M(\xi)^\dag$ such that $M(\xi)M(\xi)^\dag=Id$
$\forall \xi \neq (0,\ldots,0)$, and if
$\exists \rho <1$ such that $\forall j\in {\mathbb Z}$ and
$\forall \varepsilon\in\{ 0,1\}^d\setminus \{ (0,\ldots,0) \}$
$\exists M_{\omega_{j,\varepsilon}}\in {\mathbb R}^{n\times m}$ such that
$\forall \xi\in\prod_{i=1}^d\pm[\varepsilon_i 2^{j}\pi,(\varepsilon_i+1) 2^{j}\pi]$,
$\|Id-M(\xi)M_{\omega_{j,\varepsilon}}^{\dag}\|\leq \rho$ 
then the sequence (\ref{sequ}) using the MRA decomposition with Shannon wavelets, converges.
\end{theorem}

\noindent \emph{proof:}\\
We recall that we are in the case $b=\tilde B=1$ of theorem \ref{th1} since we deal with Shannon
wavelets.\\
The partition of the support of $\widehat{\v}$ operated by the MRA decomposition
is the following:
$$
J=\{(j,\varepsilon)\in\Z\times\{ 0,1\}^{d,*} \}
$$
$$
\v=\sum_{(j,\varepsilon)\in J} \v_j^\varepsilon \quad {\rm with} \quad
{\rm supp\,} \widehat{\v}_j^\varepsilon \subset \pm [\varepsilon_i 2^{j}\pi,(\varepsilon_i+1) 2^{j}\pi]
$$
This writting is due to the fact that the case $\varepsilon_i=0$ corresponds to a scaling function $\phi_j$
for the variable $x_i$,
and $\varepsilon_i=1$ to a wavelet function $\psi_j$. Owing the fact that
${\rm supp\,} \widehat{\phi}_j\subset [-2^j\pi,2^j\pi]$ and 
${\rm supp\,} \widehat{\psi}_j\subset [-2^{j+1}\pi,-2^j\pi]\cup [2^j\pi,2^{j+1}\pi]$, we obtain
the sets indicated in the theorem.
This case is represented in figure \ref{isotropic}.\\
Then we apply theorem \ref{th1} to obtain the convergence.

\begin{remark}
The fact that $M(\xi)$ admits a pseudo-inverse $M(\xi)^\dag$ such that $M(\xi)M(\xi)^\dag=Id$ is
implied by $\exists \rho <1,~\exists M_{\omega_{\bf j}}\in {\mathbb R}^{n\times m}$
such that $\|Id-M(\xi)M_{\omega_{\bf j}}^{\dag}\|\leq \rho$.
\end{remark}

\begin{theorem}[Tensorial wavelets]
If the symbol matix $M(\xi)$ admits a pseudo-inverse $M(\xi)^\dag$ such that $M(\xi)M(\xi)^\dag=Id$
$\forall \xi \in \mathbb{R}^d\setminus \{\xi\in\mathbb{R}^d$ such that $\xi_1\ldots\xi_d=0 \}$, and
$\exists \rho <1$ such that
$\forall {\bf j}\in {\mathbb Z}^d~~\exists M_{\omega_{\bf j}}\in {\mathbb R}^{n\times m}$ such that
$\forall \xi\in\prod_{i=1}^d\pm[2^{j_i}\pi,2^{j_i+1}\pi],~~\|Id-M(\xi)M_{\omega_{\bf j}}^{\dag}\|\leq \rho$
then the method converges with the tensorial wavelet decomposition.
\end{theorem}

\noindent \emph{proof:}\\
Anew, we use theorem \ref{th1} with $J=\Z^d$ and
$$
\v=\sum_{\j\in J} \v_\j \quad {\rm with} \quad
{\rm supp\,} \widehat{\v}_\j \subset \pm [2^{j_i}\pi,2^{j_i+1}\pi]
$$
This wavelet decomposition part the frequency domain as represented in figure \ref{anisotropic}.\\

\begin{remark}
If we consider only constant matrices operating on the wavelet coefficients, the resulting operations
on the Fourier transform of the functions are symetric by reflection along all axes:
\begin{equation}
\label{symsymb}
\forall i\in\{1,\dots,d\},~ M_{\omega}(\xi_1,\dots,-\xi_i,\dots,\xi_d)=M_{\omega}(\xi_1,\dots,\xi_i,\dots,\xi_d)
\end{equation}
On the other hand, as we deal with real functions, the approximation matrix must be real.
\end{remark}

\begin{remark}
The best aproximation $M_{\omega_{\bf j}}$ of $M(\xi)$, for the inversion is given by
$$
M_{\omega_{\bf j}}=\arg \min_{\mu\in\R^{n\times m}} \sup_{\xi\in supp(\v_\j)}\| Id- M(\xi) \mu^{\dag}  \|
$$
\end{remark}

\begin{example}
As we'll see in part \ref{laplw}, the operator $\Delta^{-1}$ matches the two cases. Wavelet algorithms
converge in the MRA context and in the tensorial one.
\end{example}

\begin{example}
If we consider the 1-D symbols $p:\R\to\C,\,\xi\mapsto p(\xi)$ that are continuous, an example
that doesn't match the conditions of the theorem is given
by a symbol $p$ such that: $p(\xi)=e^{i\xi}$ for $\xi\in [\pi,2\pi]$. Then the
operator whose symbol is $p$ can't satisfy $\exists \rho<1,\,\exists \mu \in {\mathbb R}$
such that $\forall \xi\in [\pi,2\pi],~\|1-p(\xi)\mu^{-1}\|\leq \rho$.
\end{example}

\begin{example}
On the other hand, the 1-D symbols $p:\,\mathbb{R}\to \mathbb{R},~\xi\mapsto p(\xi)$ that are continuous,
verify $\forall \xi\neq 0,~p(\xi)\neq 0$ and
$$
\sup_{j\in{\mathbb Z}}
\frac{\sup_{\xi\in [2^{j}\pi,2^{j+1}\pi]}(|p(\xi)|)}{\inf_{\xi\in [2^{j}\pi,2^{j+1}\pi]}(|p(\xi)|)}
<+\infty
$$
can be approximated by a constant $\omega_j$ on each interval $\pm[2^{j}\pi,2^{j+1}\pi]$ with optimal value verifying
$$
\omega_j^{-1}=\frac{1}{2}(m_j^{-1}+M_j^{-1})
$$
with $m_j=\inf_{\xi\in [2^{j}\pi,2^{j+1}\pi]}(|p(\xi)|)$
and $M_j=\sup_{\xi\in [2^{j}\pi,2^{j+1}\pi]}(|p(\xi)|)$\\
That is the case for functions with polynomial increase, since for $\xi^\alpha$,
$\frac{(2^{j+1}\pi)^\alpha}{(2^{j}\pi)^\alpha}=2^{\alpha}$.
\end{example}

\begin{example}
In 2-D, even for real operator matrices, the approximation by constant matrices can fail.
For instance, if we consider a symbol matrix $M$ such that:
$$
M(\xi_1,\xi_2)=\left( \begin{array}{cc} \cos(\xi_1) & -\sin(\xi_1) \\ \sin(\xi_1) & \cos(\xi_1) \end{array} \right)~~~~~
\textrm{for}~~\xi\in [\pi,2\pi]^2
$$
any wavelet approximation by constant matrices $\mu\in \R^{2\times 2}$ fails: either $\| Id+\mu \|\geq 1$ or $\| Id-\mu \|\geq 1$.
\end{example}


\section{Derivation of wavelets}
\label{derivew}

P.~G.~Lemari\'e-Rieusset \cite{Lem92} showed that derivating or integrating a biorthogonal wavelet
basis provided a new wavelet basis.
It allows us to construct two different one-dimensional multiresolution analyses
of $L^2(\mathbb{R})$ related by differentiation and integration.\\
\begin{theorem}[Derivation of wavelets]\cite{Lem92}
\label{thLem92}
Let   $(V^1_j)_{j\in{\mathbb Z}}$ be a one-dimensional MRA, 
with a differentiable
scaling function $\varphi_1$, ($V_0^1=\textrm{span}\{\varphi_1(x-k),k\in{\mathbb Z}\}$), and a
wavelet $\psi_1$. There exists a second MRA $(V^0_j)_{j\in{\mathbb Z}}$ with a scaling function $\varphi_0$
($V_0^0=\textrm{span}\{\varphi_0(x-k),k\in{\mathbb Z}\}$) and a wavelet $\psi_0$ satisfying:
\begin{equation}
\label{derb}
\varphi_1'(x)=\varphi_0(x)-\varphi_0(x-1) \qquad 
\psi_1'(x)=4~\psi_0(x)
\end{equation}
Expressed with its Fourier transform this relation writes:
$$
i\xi\widehat{\psi_1}(\xi)=4\widehat{\psi_0}(\xi)
$$
The filters $(m_0,m_0^*)$ and $(m_1,m_1^*)$ attached respectively to the MRA's $(V^0_j)_{j\in{\mathbb Z}}$ and $(V^1_j)_{j\in{\mathbb Z}}$ verify:
$$
m_0(\xi)=\frac{2}{1+e^{-i\xi}}~m_1(\xi) \quad \textrm{and} \quad
m_0^*(\xi)=\frac{1+e^{i\xi}}{2}~m_1^*(\xi)
$$
\end{theorem}
If the wavelet $\psi_1$ is $C^{n}$ and has $p$ zero moments (i.e. $\widehat{\psi}_1$ is derivable $p-1$ times in a
neighborhood of $0$ and $\widehat{\psi}_1^{(k)}(0)=0$ for $0\leq k\leq p-1$)
after such an operation, 
the wavelet $\psi_0$ has regularity $C^{n-1}$ and $p+1$ zero moments.
\begin{remark}
As the Shannon wavelets are $C^{\infty}$ and have an infinite number of zero moments, they can be derivated or
integrated in order to obtain biorthogonal
wavelets satisfying the relations (\ref{derb}) of the theorem \ref{thLem92}. And we can iterate the derivation or the integration
of these wavelets as many times as we wish in order to obtain derivatives of arbitrary order: $\dots$, $\psi_{-2}$, $\psi_{-1}$,
$\dots$, $\psi_{2}$, $\dots$, with $\psi_{0}$ the original Shannon wavelet.
\end{remark}
On account of the above remark, we can introduce a new operation thanks to the wavelet decomposition of a function ${v}$.
Indeed, if we use the tensorial wavelet decomposition, we can derivate or integrate in every directions.
For instance, if we write the wavelet decomposition of $v$ with wavelets $\psi_0$ for each tensorial components
except for $i$ for which we take $\psi_1$ where $\psi_0$ and $\psi_1$ are related by the derivation
relation (\ref{derb}) as in theorem \ref{thLem92}.
$$
v(x)=\sum_{{\bf k,j}\in{\mathbb Z}^d} d_{\bf j\,k}\psi_0(2^{j_1}x_1-k_1)\dots\psi_1(2^{j_i}x_{i}-k_{i})\dots\psi_0(2^{j_d}x_d-k_d)
$$
Then if we put for $u$:
$$
u(x)=\sum_{{\bf k,j}\in{\mathbb Z}^d} 4\cdot2^{j_i}d_{\bf j\,k}\psi_0(2^{j_1}x_1-k_1)\dots\psi_0(2^{j_i}x_{i}-k_{i})\dots\psi_0(2^{j_d}x_d-k_d)
$$
We obtain:
$$
u(x)=\frac{\partial v}{\partial x_i}(x) \quad \textrm{or, in Fourier} \quad \widehat{u}(\xi)=i\xi_i\widehat{v}(\xi)
$$


\section{Constructible approximations}

Here we restict ourselves to the case when $m=n$. From the results of the previous section, it comes:

\begin{theorem}[Set of constructible operators]
The set of symbol matrices that are constructible by multiplying the wavelet coefficients by some constants depending
on the parameter $\j$ and by derivating wavelets as
in theorem \ref{thLem92} is the $\mathbb{R}$-algebra of $\xi\mapsto {\mathbb C}^{n\times n}(\xi)$ generated by the elements
$\{(\delta_{i,j})\}_{1\leq i,j\leq n}$, $\{i\xi_i I\}_{1\leq i\leq d}$ and $\{i\xi_i^{-1} I\}_{1\leq i\leq d}$,
with $(\delta_{i,j})$ denoting the matrix which is zero everywhere except at line $i$ and column $j$ where it is $1$.
\end{theorem}
This theorem enables us to diversify our wavelet approximations of differential operators and extends the result
of section \ref{CDVtens}. But it still remains rather limited since for instance,
in dimension larger than 2, we cannot reach $\Delta^{-1}$, the inverse Laplacian, with these operations.


\section{Convergence with the wavelet packets}
\label{contwavpak}

Let $A$ be an operator from $(H^{t/2}({\mathbb R}^d))^n$ to $(H^{-t/2}({\mathbb R}^d))^n$,
having a continuous symbol $M(\xi)$ almost everywhere invertible on $\R^d$ in the sens of Riemann measure (i.e.
such that for all compact sets $K$ of $\mathbb{R}^d$,
$K\cap (\det(M))^{-1}(\{ 0 \})$ the subset of $K$ where $M(\xi)$ is not invertible
has a vanishing Riemann measure), and verifies the condition (\ref{symsymb}). Then it can be approximated
by constant matrices $M_{\omega_{\bf j}}$
with Shannon wavelet packets providing an ad hoc partition of the frequency domain.

\begin{theorem}
For all linear operator $A$ satisfying the above conditions, we can numerically solve the equation $A\u=\v$ with a wavelet packet method.
i.e. $\forall \varepsilon>0$, we can find $\u_\varepsilon$ such that $\|\v-A\u_\epsilon\|<\varepsilon$, thanks to the wavelet
algorithms described in part \ref{CDVtens}.
\end{theorem}

\noindent \emph{Proof:} First we build a finite set of rectangles $\{ \omega_\j \}_{\j\in J}$ of the type
$$
\omega_\j=\prod_{i=1}^d \left[2^{j_i}\ell_i,2^{j_i}(\ell_i+1)\right]
$$
with $\j=(j_1,\dots,j_d,\ell_1,\dots,\ell_d)$, $j_i\in{\mathbb Z}$ and $k_i\in\N$, such that if $\Omega=\bigcup_{\j\in J} \omega_\j$ then
$\| \widehat{v}-\widehat{v}_{|\Omega} \|<\varepsilon/2$,
and for $\xi_\j=(2^{j_1}\ell_1,\dots,2^{j_i}\ell_i,\dots,2^{j_d}\ell_d)\in \omega_\j$,
$$
\sup_{\j} \left( \sup_{\xi\in \omega_\j}\| Id- M(\xi)M^{-1}(\xi_\j) \| \right)<1
$$
Remark that, as a consequence, $M(\xi)$ is invertible on $\Omega$.\\
On the set $\Omega$, we can apply the wavelet algotithm of part \ref{CDVtens} with Shannon wavelet packets,
and get $\u_\epsilon$ such that
$\| \widehat{\v}_{|\Omega}-\widehat{A}_{|\Omega}\widehat{\u_\epsilon}_{|\Omega} \|<\varepsilon/2$,
then we extend $\widehat{\u_\epsilon}$ to $\mathbb{R}^d$ by taking $\widehat{\u_\epsilon}_{|^c\Omega}=0$.

\begin{remark}
This approach is valid for Shannon wavelets. But in practice we would like to use other
wavelets for which the frequency partition induced
by the wavelet packets is very difficult to control.
\end{remark}

\begin{remark}
Nevertheless, this method should help us to improve the existing algorithms (see parts \ref{laplw} and \ref{leraynum}).
\end{remark}


\section{Implicit Laplacian}
\label{laplw}

First, let us begin with the elliptic operator $A=(Id-\alpha\Delta)$, $\alpha\geq 0$.
The above study allows us to state precise properties for the
{\bf wavelet iterative algorithm} solving $(Id-\alpha\Delta){\bf u}={\bf v}$ (in Navier-Stokes $\alpha=\nu\delta t$).

\ 

First, we consider a Shannon wavelet decomposition for $\u$ and $\v$:
$$
(Id-\alpha\Delta)\sum_{j\in\mathbb{Z}^d}{\bf u}_{\bf j} =\sum_{j\in\mathbb{Z}^d}{\bf v}_{\bf j}
$$
then, for each ${\bf j}\in\mathbb{Z}^d$, we solve
$$
(Id-\alpha\Delta){\bf u}_{\bf j}={\bf v}_{\bf j}
$$
with ${\rm supp\,}(\widehat{{\bf u}_{\bf j}})=\,{\rm supp\,}(\widehat{{\bf v}_{\bf j}})\subset \prod_{i=1}^d \pm[2^{j_i+1}\pi,2^{j_i}\pi]$\\
We approximate $M(\xi)=(1+\alpha|\xi|^2)\,Id$ with $\xi\in\prod_{i=1}^d \pm[2^{j_i+1}\pi,2^{j_i}\pi]$ by
$M_{\omega_{\bf j}}=(1+\alpha\omega_{\bf j}^2)\,Id$ with $\omega_{\bf j}\in{\mathbb R}$ properly chosen.

Now, we'll see what is the most convenient value for the parameter $\omega_{\bf j}$.
The symbol matrix $(Id-M(\xi)M_{\omega_{\bf j}}^{-1})$ of the operator $(Id-A M_{\omega_{\bf j}}^{-1})$
is a diagonal matrix $\lambda Id$ whose element $\lambda\in\R$ is equal to:
$$
\lambda=1-\frac{1+\alpha|\xi|^2}{1+\alpha\omega_{\bf j}^2}
$$
This expression has a minimal maximum for $|\xi|^2\in[a_{\bf j}^2,b_{\bf j}^2]$ at $\omega_{\bf j}^2=\frac{a_{\bf j}^2+b_{\bf j}^2}{2}$.\\
If we take this value for $\omega_{\bf j}$, then
$$
\forall\xi\in\mathbb{R}^d~/~a_{\bf j}^2\leq |\xi|^2 \leq b_{\bf j}^2,~~~
|\lambda(Id-M(\xi)M_{\omega_{\bf j}}^{-1})|\leq \frac{\alpha(b_{\bf j}^2-a_{\bf j}^2)}{2+\alpha(b_{\bf j}^2+a_{\bf j}^2)}
$$

Then in the tensorial case, with the Shannon wavelets, as for fixed $\j$, $\forall i\in\{ 1,\dots,d \},~\xi_i\in \pm[2^{j}\pi,2^{j+1}\pi]$,
then $b_{\bf j}=2a_{\bf j}$ and $|\lambda|\leq \frac{3\alpha}{2a_{\bf j}^{-2}+5\alpha}$. The worst case appears for $\alpha\to +\infty$,
then the convergence rate tends to $\rho=\frac{3}{5}=0.6$.\\
In the case of the MRA, we have always $b_{i,\bf j}^2\leq (d+1)a_{i,\bf j}^2$, that allows the algorithm to converge too:
$\rho=\frac{d\alpha}{2a_{\bf j}^{-2}+(d+2)\alpha}\to_{\alpha\to +\infty}\frac{d}{d+2}$.

If we now consider what should be obtained with the Shannon wavelet packets, we take $b=3/2\,a$ in the tensorial case, then
if $\alpha\to \infty$, $\rho=\frac{5}{13}$ ($\sim 0.38$). Roughly speaking, we improve the convergence by a factor $2$ for each
wavelet packet refinement.

\ 

\noindent {\bf Observed convergence for the implicit Laplacian}\\
For the operator $Id-\alpha \Delta$, the convergence of the wavelet algorithm is very fast if $\alpha$ is small compared
to the smallest computed scale. Anyway, for the Laplacian operator $\Delta$ ($\alpha\to +\infty$) the same algorithm still converges,
and the observed convergence rate is around $0.75$ for spline wavelets of order $4$.


\section{Leray projector}
\label{leraynum}

{\bf Principle of the Helmholtz decomposition:}\\
Let $\u\in (L^2(\mathbb{R}^n))^n$ be a vector field. We can decompose $\u$ as follows:
$$
\u=\u_{\,{\rm div}} +\u_{\,{\rm curl}} ~~~{\rm where}~~~ \u_{\,{\rm div}}=\mathbf{curl}~ \psi~,~~ \u_{\,{\rm curl}}=\nabla p
$$
The functions $\mathbf{curl} ~\psi$ and $\nabla p$ are orthogonal in $(L^2(\mathbb{R}^n))^n$ and are unique.
$$
(L^2(\mathbb{R}^n))^n={\textbf{H}}_{\,{\rm div}\,0}(\mathbb{R}^n)\oplus^{\bot} {\textbf{H}}_{\,{\rm curl}\,0}(\mathbb{R}^n)
$$

In Navier-Stokes, this decomposition is very important to project the term $\u.\nabla \u$ onto
${\textbf{H}}_{\,{\rm div}\,0}(\mathbb{R}^n)$ the space of divergence free vector functions (see \cite{CM93}).

\ 

The wavelet algorithm that will be studied here was originally designed in
\cite{DP06a}. The proof of its convergence for the Shannon wavelet is new.
The {\bf wavelet iterative algorithm for ${\mathbb P}$}, the orthogonal projector (in $L^2$) on ${\textbf{H}}_{\,{\rm div}\,0}(\mathbb{R}^n)$,
uses the result on the wavelet derivation of part \ref{derivew}.
We approximate the Leray projector ${\mathbb P}$: ${\mathbb P}\widehat{\bf u}=M(\xi)\widehat{\bf u}$ with
$$
M(\xi)= Id - \frac{1}{|\xi|^2}\left[ \begin{array}{c} \xi_1  \\ \vdots \\  \xi_n \end{array} \right]
\times \left[ \xi_1 \dots \xi_n \right]
$$
by
$$
M_{\omega}=
\left( Id- \frac{1}{|\omega|^2}\left[ \begin{array}{c} \frac{\omega_1^2}{\xi_1} \\ \vdots \\ \frac{\omega_n^2}{\xi_n} \end{array} \right]
\times \left[ \xi_1 \dots \xi_n \right]  \right)\left(Id-\frac{1}{|\omega|^2}\left[ \begin{array}{c} \xi_1  \\ \vdots \\  \xi_n \end{array} \right]
\times \left[ \frac{\omega_1^2}{\xi_1} \dots \frac{\omega_n^2}{\xi_n} \right] \right)
$$
where we used the notation $\omega_i=2^{j_i}$.
\begin{remark}
One can notice that this projector $M_{\omega}$ actually projects the vector field $\u$ on the space of divergence free vector functions
(indeed $<[\xi_1 \dots \xi_n],M_{\omega}\widehat{u}(\xi)>=0$).
\end{remark}

We also have to extract the gradient part of $\u$ by the approximation $L_\omega$ of $\nabla \Delta^{-1} (\,{\rm div} \cdot)$:
$$
L_\omega=\frac{1}{|\omega|^2}\left[ \begin{array}{c} \xi_1  \\ \vdots \\  \xi_n \end{array} \right]
\times \left[ \frac{\omega_1^2}{\xi_1} \dots \frac{\omega_n^2}{\xi_n} \right]
$$
\begin{theorem}
The sequence (\ref{sequw}) writes:
\begin{equation}
\label{lerit}
{\bf v}_0=\u \quad {\rm and} \quad {\bf v}_{n+1}={\bf v}_{n}-M_{\omega}{\bf v}_{n}-L_{\omega}{\bf v}_{n}
\end{equation}
This sequence $({\bf v}_n)$ goes to $0$ in $L^2$ norm in the case of Shannon wavelets. The iteration of the operation (\ref{lerit})
allows us to separate the divergence free part $\u_{\rm div}=\sum_{n\in\N}M_\omega\,{\bf v}_{n}$
of the function $\u$ from its gradient part $\u_{\rm curl}=\sum_{n\in\N}L_\omega\,{\bf v}_{n}$.
\end{theorem}
\emph{proof:}\\
The matrix $(Id-M_{\omega}-L_{\omega})$ can be written, for $\xi\in \prod_i \pm[2^{j_i}\pi,2^{j_i+1}\pi]$
$$
(Id-M_{\omega}-L_{\omega})=\frac{1}{|\omega|^2}\left[ \begin{array}{c} \frac{\omega_1^2}{\xi_1} \\ \vdots \\ \frac{\omega_n^2}{\xi_n}
\end{array} \right] \times \left[ \xi_1 \dots \xi_n \right]\, \times \, \left( Id-L_\omega \right)
$$
All its eigenvalues are equal to zero except one that equals
$$
\lambda=1-\left(\sum_{k=1}^n \xi_k^2\right)\left(\sum_{k=1}^n \frac{\omega_k^4}{|\omega|^4\xi_k^2}\right)
$$
i.e., with $\zeta_k=\frac{\xi_k}{\omega_k}$,
$$
\lambda=1-\left(\sum_{k=1}^n \frac{\omega_k^2}{|\omega|^2}\zeta_k^2\right)\left(\sum_{k=1}^n \frac{\omega_k^2}{|\omega|^2}\zeta_k^{-2}\right)
$$
Then the Kantorovitch inequality yields
$$
|\lambda| \leq\frac{1}{4}\left(\frac{\min_k|\zeta_k|}{\max_k|\zeta_k|}+\frac{\max_k|\zeta_k|}{\min_k|\zeta_k|}\right)^2-1
$$

\begin{figure}[t]
\begin{center}
\includegraphics[scale=0.8,angle=00]{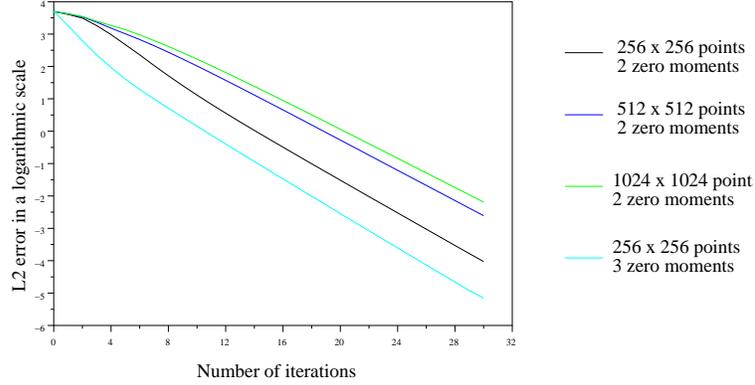}
\end{center}
\caption{\label{convwLer} Observed convergence rates for various wavelets and numbers of grid points}
\end{figure}

If $|\zeta_k|\in[a_k,b_k]$ for each $k$,
we put $a=\min_k~a_k$ and $b=\max_k~b_k$. Then
$$
|\lambda|\leq \frac{1}{4}\left(\frac{a}{b}+\frac{b}{a}\right)^2-1
$$
For the Shannon wavelets, $b=2a$, then convergence is assured since $|\lambda|\leq \frac{9}{16}$ ($\sim 0.56$).\\
In the case of the Shannon wavelet packets, $b=\frac{3}{2}a$, the convergence rate should be $\rho=|\lambda|\leq \frac{25}{144}$ ($\sim 0.17$).

\ 

\noindent {\bf Observed convergence for the Leray projector}\\
The convergence has been tested successfully on variate 2D and 3D fields.
The observed convergence rates with spline wavelets of order 2 and 3 are around 0.5 (see figure \ref{convwLer}).
The reference \cite{DP06a} provides technical explanations for the impementation of this algorithm. 


\section{Conclusion, perspectives}
  
This work provides an original point of view on the wavelet algorithms and links two achievements of
the wavelet theory:
\begin{itemize}
\item the wavelet approximation of operators like it is presented in A.~Cohen's, W.~Dahmen's and R.~DeVore's work \cite{C00,CDV01,CDV02,CHR02},
\item the divergence-free wavelet transform derived from P.-G.~Lemari\'e's and K.~Urban's works \cite{Lem92,U00}
as well as the author's \cite{Der06,DP06a,DP06b}.
\end{itemize}
The main results of this paper are the establishment of general conditions for the convergence of wavelet algorithms
with Shannon wavelets, the theoretical construction of wavelet
algorithms in order to approximate operators with constant coefficients, the exact computation of
the convergence rates and their optimisation for the implicit Laplacian operator and the Leray projector.
It also gives a simple view of the wavelet preconditionning
and a glance to what could be done thanks to wavelet packets. Indeed, wavelet packets seem to provide a powerfull solver for
many kinds of PDE's but their use is still theoretical and prospective.

The progress to achieve should go toward a more general frame for the proofs of the convergence
of wavelet algorithms that would allow us not to restrict ourselves to Shannon wavelets. As a result
we should be able to introduce wavelets on the interval in this frame.
We are still missing efficient wavelet packets concerning the frequency localisation.
An interesting perspective would be to study what happens in the case of operators with non constant coefficients

These theoretical assertions have already an application since these technics are derived to simulate the Navier-Stokes equations
(see \cite{Der06}, and a forthcoming paper in Siam Multiscale Simulation and Analysis).

\section*{Aknowledgements} 
The author thanks the University of Ulm and particularly Karsten Urban and Kai Bittner from the Numerical Analysis
team for hosting during year 2006.
He also thankfully acknowledges partial financial support from the European Union project IHP on `Breaking Complexity',
contract HPRN-CT-2002-00286.


\bibliographystyle{plain}

\begin{thebibliography}{99}
{\small


\bibitem{CM93} Chorin, A.J., and J.E.~Marsden, \emph{A Mathematical Introduction to Fluid
Mechanics}, book, 3rd ed., Springer, 1993. 


\bibitem{C00} A.~Cohen, \emph{Wavelet methods in numerical analysis},
Handbook of Numerical Analysis, vol. VII, P.G.Ciarlet and J.L.Lions eds.,
Elsevier, Amsterdam, 2000.

\bibitem{CDV01} A.~Cohen, W.~Dahmen, R.~DeVore, \emph{Adaptive Wavelet Methods for Elliptic Operator
Equations Convergence Rates}, Math. Comp. {\bf 70}, 27-75, 2001.

\bibitem{CDV02} A.~Cohen, W.~Dahmen, R.~DeVore, \emph{Adaptive Wavelet Methods for operator equations: beyond the Elliptic Case},
Found. Comput. Math., 2, no. 3 , pp. 203-245, 2002.

\bibitem{Coh03} A.~Cohen, \emph{Numerical analysis of wavelet methods}, Studies in mathematics and its applications,
Elsevier, Amsterdam, 2003.

\bibitem{CHR02} A.~Cohen, M.~Hoffmann, M.~Reis, \emph{Adaptive wavelet Galerkin methods for
linear inverse problems}, Siam J. Numer. Anal., 2002.


\bibitem{Der06} E.~Deriaz,
\emph{Ondelettes pour la Simulation des \'Ecoulements Fluides Incompressibles en Turbulence} (in french),
Th\`ese de doctorat de l'INP Grenoble, 2006.

\bibitem{DP06a}
		E.~Deriaz and V.~Perrier,
		\emph{Divergence-free Wavelets in 2D and 3D, application to the Navier-Stokes equations},
		J. of Turbulence, {\bf 7}(3): 1--37, 2006.
		
\bibitem{DP06b}
		E.~Deriaz, K.~Bittner and V.~Perrier,
		\emph{D\'ecomposition de Helmholtz par ondelettes~: convergence d'un algorithme it\'eratif},
		(in french), ESAIM : Proc, submitted, 2006.

\bibitem{Don94} D.~Donoho, \emph{De-Noising via Soft Thresholding},
IEEE Trans. Inf. Theory, {\bf 41}(3) 613-627, 1994.





\bibitem{Hor85} L.~H\"ormander, \emph{The Analysis of Linear Partial Differential Operators},
t. III, Springer Verlag, 1985.

\bibitem{Jaf92} S.~Jaffard, \emph{Wavelets methods for fast resolution of elliptic problems},
SIAM J. Numer. Anal. {\bf 29}, 965-986, 1992.



\bibitem{KL98} J.-P.~Kahane and P.G.~Lemari\'e-Rieusset, \emph{Fourier
series and wavelets}, book, Gordon \& Breach, London, 1995.

\bibitem{Lem92} P.G.~Lemari\'e-Rieusset, \emph{Analyses
multi-r\'esolutions non orthogonales, commutation entre projecteurs et
d\'erivation et ondelettes vecteurs \`a divergence nulle} (in french), Revista Matem\'atica
Iberoamericana, {\bf 8}(2): 221-236, 1992.




\bibitem{M98} S.~Mallat, \emph{A Wavelet Tour of Signal Processing}, 
Academic Press, 1998.

\bibitem{SKF97} K.~Schneider, N.~Kevlahan and M.~Farge,
\emph{Comparison of an adaptive wavelet method and nonlinearly filtered pseudo-spectral methods for two-dimensional turbulence}, 
Theor. Comput. Fluid Dyn. {\bf 9}: 191-206, 1997.

\bibitem{U00} K.~Urban, \emph{Wavelet Bases in H(div) and H(curl)}, Mathematics
of Computation {\bf 70}(234): 739-766, 2000.


}

\end{thebibliography}

\end{document}